\documentclass[11pt]{article}

\usepackage{amsmath,amssymb,amsthm,amsfonts,mathrsfs, mathtools, bbm, dsfont}
\usepackage{graphicx}

\usepackage{subfigure}
\usepackage{ifthen}
\usepackage{multirow, multicol}
\usepackage{float}
\usepackage{subeqnarray, array}
\usepackage{enumitem,booktabs}
\usepackage{xcolor}
\usepackage[font=small]{caption}
\graphicspath{{Figures/}}

\usepackage{fullpage}

\newcommand{\beq}{\begin{equation}}
\newcommand{\eeq}{\end{equation}}
\newcommand{\beqa}{\begin{eqnarray}}
\newcommand{\eeqa}{\end{eqnarray}}
\newcommand{\beqan}{\begin{eqnarray*}}
\newcommand{\eeqan}{\end{eqnarray*}}

\newcommand{\E}{\mathbb{E} }

\newcommand{\prob}{\mathbb{P}}
\newcommand{\ind}[1]{\mathbbm{1}_{\lb#1\rb}}

\newcommand{\Gset}{\mathbb{G}}

\newcommand{\Rset}{\mathbb{R}}

\newcommand{\Acal}{{\cal A}}

\newcommand{\Ecal}{{\cal E}}
\newcommand{\Fcal}{{\cal F}}
\newcommand{\Gcal}{{\cal G}}

\newcommand{\Mcal}{{\cal M}}

\newcommand{\Pcal}{{\cal P}}

\newcommand{\Gfk}{{\mathfrak{G}}}

\newcommand{\Nfk}{{\mathfrak{N}}}

\newcommand{\Rfk}{{\mathfrak{R}}}

\newcommand{\bone}{\mathds{1}}

\renewcommand{\v}[1]{{\mathbf{#1}}}

\renewcommand{\[}{\left[}
\renewcommand{\]}{\right]}
\newcommand{\lb}{\left\{}
\newcommand{\rb}{\right\}}
\renewcommand{\(}{\left(}
\renewcommand{\)}{\right)}

\newcounter{l1}
\newcounter{l2}
\newcounter{l3}
\newcommand{\bdotlist}{\begin{list}{$\bullet$}{}}
\newcommand{\bboxlist}{\begin{list}{$\Box$}{}}
\newcommand{\bbboxlist}{\begin{list}{\raisebox{.005in}{{\tiny
$\blacksquare$ \ \ }}}{}}
\newcommand{\bdashlist}{\begin{list}{$-$}{} }
\newcommand{\blist}{\begin{list}{}{} }
\newcommand{\barablist}{\begin{list}{\arabic{l1}}{\usecounter{l1}}}
\newcommand{\balphlist}{\begin{list}{(\alph{l2})}{\usecounter{l2}}}
\newcommand{\bAlphlist}{\begin{list}{\Alph{l2}.}{\usecounter{l2}}}
\newcommand{\bdiamlist}{\begin{list}{$\diamond$}{}}
\newcommand{\bromalist}{\begin{list}{(\roman{l3})}{\usecounter{l3}}}

\newtheorem{proposition}{Proposition}

\renewcommand{\Gset}{{\mathfrak{G}}}

\begin{document}

\title{\textbf{Centralized Volatility Reduction for Electricity Markets}\\
{\it{{\small{\color{red}{To appear in the International Journal of Electrical Power and Energy Systems}}}}}}

\author{Khaled Alshehri \qquad Subhonmesh Bose \qquad Tamer Ba\c sar \thanks{K. Alshehri is with the Systems Engineering Department, King Fahd University of Petroleum and Minerals, Dhahran, Saudi Arabia. Email: kalshehri@kfupm.edu.sa. S. Bose and T. Ba\c{s}ar are with the Department of Electrical and Computer Engineering, University of Illinois at Urbana-Champaign. Emails: \{boses, basar1\}@illinois.edu.}
}
 \date{} 
\maketitle

\begin{abstract}
Increased penetration of wind energy will make electricity market prices more volatile. As a result, market participants will bear increased financial risks, which impact investment decisions and in turn, makes it harder to achieve sustainable energy goals. As a remedy, in this paper, we propose an insurance market that complements any wholesale market design. Our mechanism can be run by any suitable financial entity such as an independent system operator, with the aim of reducing the financial impacts of volatile prices. We provide theoretical guarantees, analytically characterize the outcomes over a copperplate power system example, and numerically explore the same for a modified IEEE 14-bus test system.
\end{abstract}
\maketitle
\setcounter{page}{1}

\section{Introduction}
Wind energy is uncertain (difficult to forecast), intermittent (shows large ramps), and largely uncontrollable (output largely cannot be altered on command). It fundamentally differs from dispatchable generation that ``can be controlled by the system operator and can be turned on and off based primarily on their economic attractiveness at every point in time'' \cite{joskow2011comparing}. 
It has been a widely recognized fact that escalated penetration of wind will dampen electricity prices. Wind is a (near) zero marginal cost resource, and hence, alters the merit-order at the base of the stack. ``Free'' wind shifts the market supply curve to the right, leading to price reduction. Empirical evidence corroborates that hypothesis, such as the analyses by Ketterer in \cite{ketterer2014impact} for the German market, Munksgaard and Morthorst in \cite{munksgaard2008wind} for the Danish market, de Miera et al. in  \cite{de2008analysing} for the Spanish electricity market, among others. Green's model-oriented analysis \cite{green2010market} for the British market resonates the same sentiments.

A perhaps less studied effect of large-scale wind integration is its contribution to price volatility. Dispatchable (and often marginal) generators need to compensate for variations in wind availability, leading to variations in energy prices. Data from various markets support that conclusion, such as the studies by Woo et al. in \cite{woo2011impact} for ERCOT, Martinez-Anido in \cite{martinez2016impact} for New England, J{\'o}nsson et al. in \cite{jonsson2010market} for the Danish market, and Ketterer in \cite{ketterer2014impact} for the German one.\footnote{Price variations differ considerably across a day; they are positively correlated with demand, as \cite{green2010market} concludes from the British market. They also exhibit seasonal variations; these variations are greater in the summer, as the Australian market analysis in \cite{higgs2015australian} reveals.} Gerasimova \cite{gerasimova2017electricity}, studying the Nord pool (Finland, Sweden, Norway, Denmark), shows that intraday price variations in parts of Finland and Sweden -- measured in terms of the expected difference in daily on-peak and off-peak prices -- have roughly doubled during the period 2008-2016 from that in 2000-2007.
Such trends are likely to persist and perhaps grow, given the rapid growth in wind penetration.

How can market participants hedge against financial risks from these price variations? Financial instruments, such as forwards, futures, swaps, and options, can help mitigate such risks; see \cite{der, swingsurvey, kluge,Biggar,Chao,gupta} for their use in electricity markets. {\color{black} In addition to hedging, options have been shown to mitigate the effects of market power in electricity markets \cite{Allaz,Anderson,Holmberg}.} The focus of the current paper is on the design of an insurance market that complements any wholesale market design. Our design is inspired by cash-settled call options that are bilateral financial instruments to hedge volatility. 

In this paper, we propose a \emph{centralized market mechanism for insurance contracts}, where a market maker facilitates the trade of insurance contracts by `matching' buyers and sellers. Our market design is different from a traditional exchange such as the European Energy Exchange and the Chicago Board Options Exchange. Here, we allocate the collection of contracts bought among sellers with the goal to reduce the aggregate volatilities in the profits received by electricity market participants. Such a mechanism will aid electricity markets with high penetration of wind by allowing the market participants to mitigate their financial risks. The market we propose is an add-on to run in parallel with any electricity market design, and hence, does not advocate any alteration to existing dispatch and pricing of electricity markets. Our contribution complements the financial risk exchange between wind power producers proposed in \cite{shin}, but it is more general in the sense that we allow any electricity market participant to buy or sell the contracts we study. 

We propose our centralized insurance clearing mechanism in Section \ref{sec:centralized} that is run by a market maker. For our mechanism, we prove that the aggregate volatilities cannot increase, and the expected merchandising surplus remains zero even if the market maker is profit-motivated. Next, in Section \ref{sec:model}, we present a dispatch and pricing model for a two-period electricity market that we apply our insurance market to, provide conditions to guarantee strict volatility reduction for market participants.
Then, Section \ref{sec:central.example} analytically illustrates volatility reductions through a stylized copperplate power system example and demonstrates how our mechanism generalizes bilateral trading of call options. We conduct numerical experiments on the IEEE 14-bus test system  \cite{IEEE14Bus} in Section \ref{sec:central.example2} to further explain the properties of our mechanism. The paper concludes in Section \ref{sec:conclusion}. All proofs are included in the Appendix.

\subsubsection*{Notation}
\label{sec:notation}
We let $\Rset$ denote the set of real numbers, and $\Rset_+$ (resp. $\Rset_{++}$) denote the set of nonnegative (resp. positive) numbers. For $z\in\Rset$, we let $z^+ := \max\{z, 0\}$. For a random variable $Z$, we denote its expectation by $\E [Z]$, its variance by ${\sf var}[Z]$, and its cross-covariance with another random variable $X$ by ${\sf cov}(X,Z)$; note that ${\sf cov}(X,X) = {\sf var}[X]$. For an event $\Ecal$, we denote its probability by $\prob\{\Ecal\}$ for a suitably defined probability measure $\prob$. The indicator function for an event $\Ecal$ is given by $\ind{\Ecal}$.
In any optimization problem, a decision variable $x$ at optimality is denoted by $x^*$.
%

\section{Centralized clearing of insurance contracts}
 \label{sec:centralized} 
 
 Consider a wholesale electricity market with multiple consumers and producers. The consumers are utility companies or retail aggregators who represent a collection of retail customers. In this model, we consider two types of producers -- dispatchable generators and variable renewable wind power producers. Dispatchable generators can alter their power output within their capabilities on command, e.g., nuclear, coal, natural gas, biomass or hydro power based power plants. In contrast, the available production capacity of variable producers rely on an intermittent resource like wind energy. The SO implements a centralized market mechanism to balance demand and supply of power within the network constraints. 

  To motivate the design of our insurance market,  consider a two-stage electricity market model. Identify $t=0$ as the forward stage, prior to the uncertainty being realized, and  $t=1$, the real-time stage. Let $(\Omega, \Fcal, \prob)$ denote the probability space describing the uncertainty. Here, $\Omega$ is the collection of possible scenarios at $t=1$,
 $\Fcal$ is a suitable $\sigma$-algebra over $\Omega$, and $\prob$ is a probability distribution over $\Omega$.
We assume that $\Omega$ is compact, and that all market participants know $\prob$.

An insurance contract in our context allows its buyer the right to claim a monetary reward equal to the positive difference between the real-time price $p^{\omega,*}$ and the strike price $K$ of an underlying commodity for an upfront fee. Adopting a game-theoretic framework, consider the case where  player $r$ is the buyer and another player $g$ is the seller. The contract costs $r$ a fee of $q\Delta$, where $q$ is the upfront price and $\Delta$ is the quantity. Once they agree on the trade triple $\(q,K,\Delta\)$, their profits in scenario $\omega$ are respectively given by
\begin{align}
\begin{aligned}
\Pi^\omega_r(q, K, \Delta) &:= \pi_r^\omega - q \Delta + \left( p^{\omega,*} - K \right)^+\Delta,\\
\Pi^\omega_g(q, K, \Delta) &:= \pi_g^\omega + q\Delta - \left( p^{\omega,*} - K \right)^+\Delta.
\end{aligned}
\label{eq:Pi}
\end{align}
In each expression, the first term is the profit from the electricity market, and the other two terms capture the aggregate return from the insurance contract. Such contracts allow market participants to reduce their profit volatilies, which are measured here in terms of their variances. That is, with a well-designed contract, one would expect ${\sf var}[\Pi^\omega_i] \leq {\sf var}[\pi^\omega_i]$ for market participant $i$. 

Peaker power plants are not always asked to produce in real time, but they are critical for resource adequacy in wholesale markets \cite{ramp3}. The California Independent System Operator (CAISO) and the Midcontinent ISO (MISO) have proposed flexible ramping products and buy capacities from peaker power plants to provide them incentives to remain online. The above contracts can provide financial incentives for peaker power plants to stay in the market, without requiring the system operator to purchase such capacities. For example, a peaker power plant $g$ can participate as a seller, and receive profit of $q\Delta$ in the forward stage. On the other hand, wind power producers face the risks of not being able to produce much power in real time, and hence, they can become buyers and guarantee a reward in such events. This leads to incentives for both wind producers and peaker plants to engage in such trades. 

Contracts of the form in \eqref{eq:Pi} are often traded bilaterally between market participants in the form of cash-settled call options, but in a wholesale market with a collection of dispatchable generators $\Gfk$ and variable generators $\Rfk$, one can conceive of $\lvert \Gfk \rvert \cdot \lvert \Rfk \rvert $ bilateral trades. It is difficult to convene and settle a large number of bilateral trades on a regular basis. Financial exchanges provide an alternative that typically seek to maximize the surplus from trading (options and other financial derivatives) with a collection of market participants, without explicit consideration of aggregate volatilities.
In this paper, we take an alternate view, and propose a centralized clearing mechanism for both buyers and sellers of insurance contracts with the goal to reduce profit volatilities in electricity markets. Such an approach leads to critical outcomes: it makes volatility reduction accessible to any market participant, does not alter dispatch and pricing of the wholesale market, and is therefore compatible with any existing electricity market design.

Consider a market maker $\Mcal$ who acts as an aggregate buyer for a collection of sellers $\Gfk$, and acts as a seller for the  buyers $\Rfk$. The SO or a suitable financial institution can fulfill the role of such an intermediary. In this paper, we primarily focus on the case where $\Mcal$ is social (e.g., $\Mcal$ is the SO). We later discuss how the problem would change if $\Mcal$ is profit-motivated.
We now describe the step-by-step procedure for clearing the insurance market by a social $\Mcal$. 
\noindent\underline{Forward stage:}

\begin{itemize}[leftmargin=3mm]

\item $\Mcal$ broadcasts a set of allowable trades $\Acal_0$, given by
$$\Acal_0 := \[0, \overline{q}\] \times \[0, \overline{K}\] \times \[0, \overline{\Delta}\] \subset \Rset^3_{+},$$
to all market participants $\Gfk \cup \Rfk$. 

\item Each $i \in \Gfk \cup \Rfk$ submits an acceptable (compact) set of trades, denoted by $\Acal_i \subseteq \Acal_0$.\footnote{$\Mcal$ can fix a parametric description of $\Acal$'s, and market participants report their parameter choices.}

\item $\Mcal$ correctly conjectures the real-time prices $p^{\omega,*}$ in each scenario $\omega$. {\color{black} Also, $\Mcal$ knows the profit functions $\pi^\omega_i$'s of all market participants for each scenario $\omega\in\Omega$. \footnote{One can consider revenues instead of profits in the design of this market. If $\Mcal$ is the SO, then the revenues are known exactly, while the profits are only known approximately. We sidestep such nuances and focus on the profit-based market throughout.}} Let the merchandising surplus (sum of profits) for $\Mcal$ be denoted by ${\sf MS}^\omega$ in scenario $\omega$. As an aggregate buyer and seller, it is given by
\begin{align*}
{\sf MS}^\omega &= \sum_{r \in \Rfk} q_r \Delta_r - \sum_{g \in \Gfk} q_g \Delta_g -\sum_{r \in \Rfk} \left(p^{\omega,*}_r - K_r\right)^+ \Delta_r + \sum_{g \in \Gfk} \left(p^{\omega, *}_g - K_g\right)^+ \delta_g^{\omega}.
\end{align*}
The market maker $\Mcal$ solves the following stochastic optimization problem to clear the insurance market.%
\begin{align}
\label{eq:DA.opt}
\begin{alignedat}{5}
&{\text{minimize}}   \ \ \  \sum_{i \in \Gfk \cup \Rfk}{\sf var}[\Pi^\omega_i], \\
& \text{subject to} \\
& \qquad  \sum_{g \in \Gfk} \Delta_g = \sum_{r \in \Rfk} \Delta_r, \\
& \qquad (q_g, K_g, \Delta_g)\in\Acal_g, \  (q_r, K_r, \Delta_r)\in\Acal_r,\\
& \qquad \hspace*{-12pt}\begin{rcases}
&\delta^\omega_g \in [0, \Delta_g], \\
&\displaystyle\sum_{g\in\Gfk} \delta^\omega_g =\sum_{r \in \Rfk} \Delta_r \ind{p^{\omega,*}_r\geq K_r},\\
&\displaystyle {\sf MS}^\omega=0,\end{rcases} \quad \prob-\text{a.s.},\\
& \qquad\text{for each } \ g \in \Gfk, \ r \in \Rfk,
\end{alignedat}
\end{align}
{\color{black} over $\(q_r,  K_r, \Delta_r \) \in \Rset^3_+$ for each $r \in \Rfk$, $\(q_g, K_g, \Delta_g \) \in \Rset^3_+$ and the $\Fcal$-measurable square-integrable maps $\delta^\omega_g: \Omega \rightarrow [0, \Delta_g]$ for each $g \in \Gfk$. 
} Here, $p^{\omega,*}_r$ denotes the market price faced by $r$. We define the same for $g$, accordingly.

{\color{black} The constraints in \eqref{eq:DA.opt} dictate that the volume of insurance contracts bought equals the amount that is sold, all trades are acceptable to market participants, and real-time payments cashable in each scenario can be allocated to the sellers. 
Imposing ${\sf MS}^\omega = 0$ ensures that the market maker maintains zero balance from insurance contracts, and purely facilitates the trade among the market participants, i.e., we adopt the viewpoint of a social $\Mcal$. The objective is to reduce profit volatilities in aggregate among acceptable trades.
} 
\item Buyer $r$ pays $q_r^* \Delta^*_r$ to $\Mcal$.
\item $\Mcal$ pays $q_g^* \Delta^*_g$ to seller $g$.
\end{itemize}

\noindent\underline{Real-time stage:}

\begin{itemize}[leftmargin=3mm]
\item Scenario $\omega$ is realized, and the real-time electricity prices $p^{\omega,*}$ become known. 

\item $\Mcal$ pays $\left(p^{\omega,*}_r - K^*_r\right)^+ \Delta^*_r$ to buyer $r$.

\item Seller $g$ pays $\left(p^{\omega, *}_g - K^*_g\right)^+\delta_g^{\omega, *}$ to $\Mcal$.

\end{itemize}

Our mechanism is guaranteed not to increase (and possibly decrease) the aggregate volatility of profits, i.e., problem \eqref{eq:DA.opt} admits an optimal solution that satisfies
\begin{equation}  \sum_{i \in \Gfk \cup \Rfk}{\sf var}[\Pi^{\omega,*}_i] -  \sum_{i \in \Gfk \cup \Rfk}{\sf var}[\pi^\omega_i] \leq 0. \label{eq:aggregate} \end{equation}
{\color{black} The above inequality follows from the properties of problem \eqref{eq:DA.opt}. The constraint set is compact and the objective function is continuous, and hence, 
 Weierstrass Theorem (see \cite{Luenberger}) guarantees the existence of an optimum. Note that $\Delta$'s being zero is always a feasible choice at which \eqref{eq:aggregate} is met with an equality.
Thus, at the optimal solution, aggregate volatility can only be lower than that at the feasible point with no trades.}

We remark that electricity market-specific considerations can be incorporated in our design if $\Mcal$ is the SO. For example, the strike prices and fees for buying/selling our contracts can be deemed to be nodally uniform, i.e., they are required to satisfy $q=q_r=q_g$ and $K=K_r=K_g$ for all market participants at a particular bus in the power system.

\subsection{How participant $i$ decides $\Acal_i$}
Consider a seller $g \in \Gfk$ who expects a profit $\pi_g^\omega$ in scenario $\omega$. From the electricity and insurance market, she receives a payoff of 
$$ \pi_g^\omega + q_g \Delta_g - \( p^{\omega,*}_g - K_g \)^+ \delta_g^\omega$$
in scenario $\omega$ with the trade triple $\(q_g, K_g, \Delta_g \)$, if $\Mcal$ allocates $\delta_g^\omega \in [0, \Delta_g]$. Having no control over $\delta_g^\omega$, assume that $g$ conjectures the worst case outcome $\delta_g^\omega = \Delta_g$ that minimizes her payoff, given by
$$\Pi^\omega_g \( q_g, K_g, \Delta_g \) := \pi_g^\omega + q_g \Delta_g - \( p^{\omega,*}_g - K_g \)^+ \Delta_g.$$

Evidence from electricity markets suggests that participants are often risk averse, e.g., see \cite{Bjorgan2,der}. To illustrate how the acceptability sets can be defined for risk-averse players, assume that a market participant perceives risk via the \emph{conditional value at risk} functional (see \cite{cvar,risk}), and finds a trade triple $\( q, K, \Delta \)$ acceptable, if 
\begin{align}
{\sf CVaR}_{\alpha_i} \[- \Pi^\omega_i \(q, K, \Delta \)\] \leq {\sf CVaR}_{\alpha_i} \[ -\pi^\omega_i \],
\label{eq:cvar.i}
\end{align}
where  $\pi^\omega_i$ and $\Pi^\omega_i$ describe her profits from the energy market and the energy-cum-insurance market, respectively, in scenario $\omega$. The {\sf CVaR} risk measure is given by
$$ {\sf CVaR}_\alpha \[z^\omega\] := \min_{u\in\Rset} \lb u + \frac{1}{1-\alpha} \E\[ \(z^\omega - u \)^+\]\rb$$
for an $\Fcal$-measurable map $z$. Parameter $\alpha \in [0,1)$ encodes the extent of risk-aversion.  If $z^\omega$ is the monetary loss in scenario $\omega$ with a smooth distribution, then ${\sf CVaR}_\alpha \[z^\omega \]$ equals the expected loss over the $1-\alpha$ fraction of the scenarios that result in the highest losses.   

If all market participants pick $\alpha=0$, it follows that 
$$   \sum_{i \in \Gfk \cup \Rfk} \E \[\Pi^\omega_i \( q_i, K_i, \Delta_i \) \] - \sum_{i \in \Gfk \cup \Rfk} \E \[ \pi^\omega_i \] \geq 0,$$
which is equivalent to the expected merchandising surplus being nonpositive ($\E\[{\sf MS}^\omega\]\leq0$). The constraints in \eqref{eq:DA.opt}, however, impose $\E\[{\sf MS}^\omega\]=0$, and hence it follows that  $$ \E \[\Pi^\omega_i \( q_i, K_i, \Delta_i \) \] = \E \[ \pi^\omega_i \]$$ for each $i\in\Gfk \cup \Rfk$. {\color{black} Our insurance market design is not limited to the above description of risk-preferences. Market participants can freely choose the trades they find acceptable via $\Acal_g$'s and $\Acal_r$'s. We illustrate the effects of risk-aversion later in Section \ref{sec:central.example}.} 
 {\color{black} Our mechanism is centralized with complete information, i.e., the market maker needs to know profit functions, risk-preferences, and correct price conjectures. While this might be challenging to implement in practice as is, our theoretical guarantees, stylized examples, and numerical experiments demonstrate that it can achieve significant volatility reductions, and hence, it may serve as an optimal benchmark that can be used to derive insightful policy recommendations for the development of new insurance markets.}

{\color{black}
\subsection{Electricity markets with multiple ex-post stages}
The proposed insurance market can run in parallel with wholesale electricity markets that have multiple ex-post stages. For example, consider an electricity market with a forward stage at $t=0$ (think day-ahead market), and multiple ex-post stages $t = 1,\ldots,T$ with $T>1$ (e.g., the multiple real-time markets). Here, $\Omega$ is the collection of possible scenarios at $t\in \{1,\dots,T\}$,
 
Denote the price for electricity faced by $i$ at $t\geq 1$ by $p^{\omega,*}_i(t)$, where $\omega$ encodes a random trajectory of available renewable supply. 
The insurance market can proceed as described, where the price signal $p^{\omega,*}_i$ is computed as the average electricity price over $T$ periods as $$p^{\omega,*}_i:=\frac{1}{T} \sum_{t=1}^T p_i^{\omega,*}(t).$$
The profit to each market participant in the objective of \eqref{eq:DA.opt} becomes the total profit over $T$ periods. That is, for a market participant $i$, denoting her profit from the electricity market at stage $t$ by $\pi^\omega_i (t)$, her total profit becomes $$\pi^\omega_i := \sum_{t=1}^T \pi^\omega_i (t).$$  
The insurance market is then defined with the parameters $\pi^\omega_i$ and $p^{\omega,*}_i$ for each $i$.}

\subsection{When the market maker is a profit-maximizer}
The insurance market mechanism in \eqref{eq:DA.opt} assumes a social intermediary. Next, consider a selfish market maker who aims at maximizing its expected merchandising surplus, and solves
\begin{align}
\label{eq:DA.opt_profit}
\begin{alignedat}{5}
&{\text{maximize}}   \ \ \  \E[{\sf MS}^\omega], \\
& \text{subject to} \\
& \qquad  \sum_{g \in \Gfk} \Delta_g = \sum_{r \in \Rfk} \Delta_r, \\
& \qquad (q_g, K_g, \Delta_g)\in\Acal_g, \  (q_r, K_r, \Delta_r)\in\Acal_r,\\
& \qquad \hspace*{-12pt}\begin{rcases}
&\delta^\omega_g \in [0, \Delta_g] \\
&\displaystyle\sum_{g\in\Gfk} \delta^\omega_g =\sum_{r \in \Rfk} \Delta_r \ind{p^{\omega,*}_r\geq K_r}\\
\end{rcases} \quad \prob-\text{a.s.},\\
& \qquad\text{for each } \ g \in \Gfk, \ r \in \Rfk.
\end{alignedat}
\end{align}
Although the market maker here is motivated to maximize profit, our next result says that a selfish intermediary does not make profits on the average! However, volatility reduction such as that in \eqref{eq:aggregate} remains challenging to derive when $\Mcal$ is a profit-maximizer.
\begin{proposition}
\label{prop:profit}
If each player $i\in\Rfk \cup \Gfk$ picks $\alpha=0$, then, $ \E[{\sf MS}^{\omega,*}]=0$ at an optimal solution of (\refeq{eq:DA.opt_profit}).
\end{proposition}


\section{Application: A Stylized Two-Stage Electricity Market Model} \label{sec:model}
In this section, we apply our mechanism to a simple, yet illustrative electricity market model to demonstrate its properties. 

\subsection{Modeling the market participants}
Consider a power network for which $\Nfk$ denotes the set of buses, and let $d^\omega_n$ denote the aggregate real-time demand in scenario $\omega$ at node $n\in\Nfk$.
Let $\Gfk$ and $\Rfk$ denote the collection of dispatchable generators and variable renewable wind power producers, respectively. {\color{black}We denote the collection of generators at node $n$ by $\Gfk_n\subseteq\Gfk$. We similarly define $\Rfk_n\subseteq\Rfk$.} We model their individual capabilities as follows. Let each dispatchable generator $g \in \Gfk$ produce $x^\omega_g$ in scenario $\omega\in\Omega$ in real time. We model its ramping capability by letting $ \lvert x_g^\omega  - x_g^{0} \rvert \leq \ell_g,$
where $x_g^{0}$ is a generator set point at the stage prior to the real-time stage, and $\ell_g$ is the ramping limit. Let the installed capacity of generator $g$ be ${x}^{\text{cap}}_g$, and hence $x_g^\omega \in [0, {x}^{\text{cap}}_g]$. Its cost of production is given by the smooth convex increasing map $c_g : [0, {x}^{\text{cap}}_g] \to \Rset_+$. Each variable renewable wind power producer $r \in \Rfk$ produces $x_r^\omega$ in scenario $\omega\in\Omega$ in real time. It has no ramping limitations, but its available production capacity is random, and we have $ x_r^\omega \in [0, \overline{x}_r^\omega] \subseteq [0, {x}^{\text{cap}}_r].$
That is, $\overline{x}_r^\omega$ denotes the random available capacity of production, and $ {x}^{\text{cap}}_r$ denotes the installed capacity for $r$. {\color{black} The cost of production for $r$ is generally linear \cite{windcost} and hence, we consider it to be a smooth convex increasing map $c_r : [0, {x}^{\text{cap}}_r] \to \Rset_+$.}
We call a vector comprised of $x_g$ for each $g \in \Gfk$ and $x_r$ for each $r\in\Rfk$ a \emph{dispatch}. The SO decides the dispatch decisions and the compensations of all market participants. We adopt the dispatch and pricing model described below, which serves as a caricature of real electricity markets \cite{Stoft,Kirschen,Morales1, Morales2}. We also adopt the commonly used DC approximation of the power flow \cite{grainger}. That is, if the supply vector is denoted by $\v{x}$, and the demand vector is denoted by $\v{d}$, then the injection $\v{x}-\v{d} \in \Pcal$, where $\Pcal$ is the injection polytope defined as
$$\Pcal := \{\v{y} \vert  \ H\v{y} \leq L, \ \bone^T \v{y} = 0\},$$
where $H$ is the shift factor matrix and $L$ denotes the vector of capacities of the transmission lines.

\subsection{{\color{black}The dispatch and the pricing model}}
Assume that SO knows $c_g, x_g^{\text{cap}}$ for each $g\in\Gfk$ and ${\color{black}c_r}, x_r^{\text{cap}}, \overline{x}_r$ for each $r \in \Rfk$, and we have the following two stages.\footnote{In practice, the cost functions are derived from supply offers.}

\subsubsection{The forward stage}

The SO computes a forward dispatch against a point forecast of all uncertain parameters. In particular, the SO replaces the random available capacity $\overline{x}^\omega_r$ by a certainty surrogate $\overline{x}^{\text{CE}}_r \in [0, \overline{x}^{\text{cap}}_r] $ for each $r \in \Rfk$. A popular surrogate\footnote{See \cite{Morales2} for an alternate certainty surrogate.} is given by $ \overline{x}^{\text{CE}}_r := \E[\overline{x}_r^\omega]$. Denote the forward dispatch by $X_g \in \Rset, \ g\in \Gfk$, and $X_r \in \Rset, \ r \in \Rfk$. This dispatch is the solution of the following optimization problem, in which the system operator minimizes the aggregate cost of production needed to meet the demand, given the forecasts of all random variables, and subject to the network constraints. 
\begin{equation*}
\begin{alignedat}{5}
&{\text{minimize}}   \ \ \ 
	& & \sum_{g \in \Gfk} c_g(X_g) + \sum_{r \in \Rfk} c_r(X_r), \\
& \text{subject to}  & & \sum_{g \in \Gfk_n} X_g + \sum_{r \in \Rfk_n} X_r = \E[d^\omega_n],\ \v{X}-\v{d}\in\Pcal, \\
&&& X_g \in [0, {x}_g^{\text{cap}}], \ \ X_r \in [0, \overline{x}^{\text{CE}}_r],\\
&&& \text{for each } g \in \Gfk, \ r \in \Rfk, \ n\in \Nfk.
\end{alignedat}
\label{eq:DAM}
\end{equation*}
The forward price at node $n$ is given by the optimal Lagrange multiplier of the energy balance constraint. Denoting this price by $P^*_n$, generator $g \in \Gfk_n$ is paid $P_n^* X_g^*$, while producer $r \in \Rfk_n$ is paid $P_n^* X_r^*$. Aggregate consumer pays $P_n^*\E[d_n^\omega]$.

\subsubsection{At real time}
Scenario $\omega$ is realized. {\color{black} Denote the real-time dispatch by $x_g^\omega \in \Rset, \ g\in \Gfk$, and $x_r^\omega \in \Rset, \ r \in \Rfk$. This dispatch is the solution of the following optimization problem, in which the system operator minimizes the aggregate real-time cost of production, subject to supply-demand balance and generation and network constraints.} 
\begin{equation*}
\begin{alignedat}{5}
&{\text{minimize}}   \ \ \ 
	& & \sum_{g \in \Gfk} c_g(x_g^\omega) + \sum_{r \in \Rfk} c_r(x_r^\omega), \\
& \text{subject to}  & & \sum_{g \in \Gfk_n} x_g^\omega + \sum_{r \in \Rfk_n} x_r^\omega = d^\omega_n,\  \v{x}^\omega-\v{d}^\omega\in\Pcal, \\
&&& x_g^\omega \in [0, {x}_g^{\text{cap}}], \lvert x_g^\omega  - X_g^* \rvert \leq \ell_g,\\
&&& x_r^\omega \in [0, \overline{x}^\omega],\ \ \text{for each } g \in \Gfk, \ r \in \Rfk, \ n\in \Nfk.
\end{alignedat}
\label{eq:RTM}
\end{equation*}
The real-time (or spot) price is again defined by the optimal Lagrange multiplier of the energy balance constraint, and is denoted by $p^{\omega,*}_n$, for each $n\in\Nfk$. Note that $X_g^*$ computed at $t=0$ defines the generator set-points $x_g^0$ for each generator $g\in\Gfk$. Generator $g \in \Gfk_n$ is paid $p^{\omega,*}_n \left( x_g^{\omega,*}- X_g^* \right)$, while producer $r \in \Rfk_n$ is paid $p^{\omega,*}_n \left( x_r^{\omega,*}- X_r^* \right)$. {\color{black} The aggregate consumer pays $p_n^{\omega,*} \left(d_n^\omega-\E[d_n^\omega] \right)$.} Demand forecasts in practice are typically quite accurate and hence, we assume $d^\omega_n=d_,$ for each $n$ for all $\omega$. {\color{black} The payments in realtime correspond to balancing energy needs in real time; the forward stage compensates for the bulk energy transactions.}

The total payments to each participant is the sum of her forward and real time payments. We denote the profits corresponding to these payments $\pi_g^\omega$ for each $g \in \Gfk$ and $\pi_r^\omega$ for each $r \in \Rfk$ in scenario $\omega$. 
The above benchmark dispatch model generally defines a suboptimal forward dispatch in that the generator set-points are  \emph{not} optimized to minimize the expected aggregate costs of production \cite{bose}. Several authors have advocated a so-called stochastic economic dispatch model, wherein the forward set-points are optimized against the expected real-time cost of balancing (cf. \cite{Wong, singleset, Bouffard1, Bouffard2}). Our insurance market design can operate in parallel to such an electricity market, and this model only serves to illustrate the properties of our mechanism.

\subsection{Strict volatility reduction for each participant}
Our mechanism reduces volatility of market participants in aggregate, but it does not guarantee that each participant reduces her volatility. Here, we provide conditions under which strict reduction in volatility is guaranteed for each participant. 
\begin{proposition}
\label{main.2}
When the centralized insurance market mechanism in \eqref{eq:DA.opt} is applied to the two-stage electricity market model, variance in profits of participant $i \in \Gfk \cup \Rfk$ reduces if and only if ${\sf cov}(2A^\omega_i+B^\omega_i,B^\omega_i) < 0$, where
\begin{align*}
&A^\omega_r=p^{\omega,*}_r(x^{\omega,*}_r-X^*_r)-c_r(x^{\omega,*}_r), \ B^\omega_r= (p^{\omega,*}_r-K_r)^+\Delta_r,\\
&A^\omega_g=p^{\omega,*}_g(x^{\omega,*}_g-X^*_g)-c_g(x^{\omega,*}_g),\ B^\omega_g= -(p^{\omega,*}_g-K_g)^+\delta^\omega_g.
\end{align*}
for each $r \in \Rfk$ and $g \in \Gfk$.
\end{proposition}
Proposition \ref{main.2} reveals that there is reduction in the volatility of a participant when the total profits in real time (from energy and insurance markets) are {\it anti-correlated} with the profits from the insurance market alone. It aligns with the intuition that variance will decrease when the insurance market supplements the profits from the energy market.


{\color{black} \section{Copperplate Power System Example}\label{sec:central.example}} 
We present here a stylized single-bus power system example (adopted from \cite{bose}) and illustrate how a bilateral trade can reduce the volatility in payments of market participants, and even mitigate the risks of financial losses for some. 

Consider a power system with two dispatchable generators and a single variable renewable wind power producer serving a demand $d$. In this example, $ \Gset := \{ B, P \}, \quad \text{and} \quad \mathfrak{R} : = \{ W \}$,
where $B$ is a base-load generator, $P$ is a peaker power plant, and $W$ is a wind power producer. 

Let $x_B^{\text{cap}} = x_P^{\text{cap}} = \infty, \quad \text{and} \quad \ell_B = 0, \ \ell_P= \infty.$
Therefore, $B$ and $P$ have unlimited generation capacities. $B$ does not have the flexibility to alter its output in real time from its forward set-point. In contrast, $P$ has no ramping limitations. For simplicity, let $B$ and $P$ have linear costs of production. $B$ has a true marginal cost $0<\epsilon<1$, and offers a unit marginal cost. $P$ has a true unit marginal cost, and offers a higher cost $1/{\rho}$, where  $\rho \in (0,1]$, i.e., generators offer higher prices than their true costs, which is an observed phenomenon in electricity markets \cite{Green}.  
Encode the uncertainty in available wind in the set
$$ \Omega := [\mu - \sqrt{3}\sigma, \mu + \sqrt{3}\sigma] \subset \Rset_+,$$
and take $\prob$ to be the uniform distribution over $\Omega$. That is, available wind is uniform with mean $\mu$ and variance $\sigma^2$. Scenario $\omega \in \Omega$ defines an available wind capacity of $\overline{x}_r^\omega = \omega$. Further, assume that $W$ produces power at zero cost, and fixed demand $d \geq \mu + \sqrt{3} \sigma$. 

This stylized example is a caricature of electricity markets with deepening penetration of variable renewable wind supply. Base-load generators, specifically nuclear power plants, have limited ramping capabilities. Natural gas based peakers can quickly ramp their power outputs. Utilizing them to balance variability can be costly. Finally, (aggregated) demand is largely inflexible but predictable. In the remainder of this section, we analyze the effect of a bilateral contract on the market outcomes for this example.

The benchmark dispatch model yields the following forward and real-time dispatch decisions $X^*$, $x^{*, \omega}$, and the forward and real-time prices $P^*, p^{*, \omega}$, respectively. See \cite{bose} for the calculations.
\begin{gather*}
X_B^* = d - \mu, \ X_P^* = 0, \ X_W^* = \mu, \\
x_B^{\omega,*} = d - \mu, \ x_W^{\omega,*} = \min\{\omega, \mu\}, \ x_P^{\omega,*} = (\mu - \omega)^+,
\\
P^* = 1, \ p^{\omega,*} = (1/\rho) \ind{\omega \in \[ \mu-\sqrt{3}\sigma, \mu \]}. 
\end{gather*}

The above dispatch and the prices yield the following profits for market participants in scenario $\omega$:
\begin{gather*}
\pi_B^\omega = (d-\mu)(1-\epsilon),  \
 \pi_P^\omega = (\mu-\omega)^+(1/\rho-1), \
\pi_W^\omega = \mu - ( \mu - \omega)^+/\rho. 
\end{gather*}

{\color{black}A keen reader would recognize that the payment to market participants from our  mechanism in \eqref{eq:DA.opt} shares parallels to that from cash-settled call options. Through the copperplate power system example, we illustrate that our market mechanism is indeed a generalization of a bilateral call option trade between $P$ and $W$. In fact, we show that our centralized mechanism is able to achieve the maximum aggregate volatility reduction among all possible bilateral call option trades between $P$ and $W$.} 

\subsection{Bilateral insurance contract between $W$ and $P$}\label{sec:motivate.bilateral}

We model a bilateral contract between $P$ and $W$ as a robust Stackelberg game (see \cite{basar}) $\Gcal$ as follows. Right after the day ahead market is settled at $t=0$, $P$ announces a premium $q \in \Rset_+$ and a strike price $K \in \Rset_+$. Then, $W$ responds by purchasing $\Delta \in [0,\sqrt{3}\sigma]$. Note that we impose  an upper bound equivalent to the maximum possible shortfall, which allows for ease of exposition. 
We say $(q^*, K^*, \Delta^*(q^*, K^*))$ constitutes a \emph{Stackelberg equilibrium}, if 
$$\E\left[\Pi_P^\omega(q^*, K^*, \Delta^*(q^*, K^*))\right] \geq \E\left[\Pi_P^\omega(q, K, \Delta^*(q, K))\right],$$
where $\Delta^*:\Rset^2_+\rightarrow[0,\sqrt{3}\sigma]$ is the best response of $W$. Given $(q,K)$, the best response $\Delta^*$ satisfies $$\E\left[\Pi_W^\omega(q, K, \Delta^*(q, K))\right] \geq \E\left[\Pi_W^\omega(q, K, \Delta(q, K))\right].$$

{\color{black} This two-player Stackelberg game has one leader and one follower, where the leader $P$ acts first and then the follower $W$ responds. $P$ chooses $(q,K)$, anticipating the best response $\Delta^*(q,K)$ by $W$ to the prices $(q,K)$ offered by $P$. Note that the response function $\Delta^*(q,K)$ might not be unique, and in this case, $P$ might consider the worst-case scenario $\Delta^*(q,K)$, and hence, it is a robust Stackelberg game \cite{basar}. We have the following result.}

\begin{proposition}
\label{prop:bilat}
The Stackelberg equilibria of $\Gcal$ are given by $(q^*,K^*) \in \Rset^2_+$ and $\Delta^*: \Rset^2_+ \to [0, \sqrt{3} \sigma]$ that satisfy one of the following two conditions:
\begin{itemize}
\item[(i)] $2q^* + K^* > \rho^{-1}$, and $\Delta^* = 0$,
\item[(ii)] $2q^* + K^* = \rho^{-1}$, and $\Delta^*\in[0, \sqrt{3} \sigma]$.
\end{itemize}
Over all equilibria with $\Delta^* = \sqrt{3} \sigma$,
 \begin{align*}
 {\sf var}\left[\Pi^\omega_W(q^*, K^*, \sqrt{3}\sigma)\right]
- {\sf var}\left[\pi^\omega_W\right] 
&= - \frac{3}{2}{{q^*K^*}}{\sigma^2}  < 0,\\
 {\sf var}\left[\Pi^\omega_P(q^*, K^*, \sqrt{3}\sigma)\right]
- {\sf var}\left[\pi^\omega_P\right] 
&= - \frac{3}{2}{{q^*(K^*-1)}}{\sigma^2}.
\end{align*}
\end{proposition}

The first kind of equilibria describes the degenerate case, where $\Delta^*=0$. $P$ and $W$ engage in trading at the second kind of equilibria, where $2q^* + K^* = \rho^{-1}$. For $\Delta^*=\sqrt{3}\sigma$, note that the bilateral trade always guarantees volatility reduction for the wind producer $W$. For $P$, reductions are only guaranteed for the equilibria that satisfy $K^*>1$. Furthermore, in expectation, profits are unchanged. Similar conclusions can be drawn for any $\Delta^*\neq0$. This stylized example illustrates how insurance contracts can help with volatility reductions, but also shows that further improvements can be attained, which can be  done by applying our centralized mechanism, which yields the largest  volatility reductions. 
\subsection{Outcomes of the centralized mechanism}
Consider an insurance market with buyer $W$ and seller $P$ and intermediary $\Mcal$. Let the price cap be given by the maximum real-time price $1/\rho$, and the trade volume be capped at $\sqrt{3} \sigma$, the maximum energy shortfall in available wind from its forward contract. Said otherwise, $\Mcal$ restricts trade triples to the set $\Acal_0 := [0, 1/\rho] \times [0, 1/\rho] \times [0, \sqrt{3}\sigma]$.

\begin{figure}
\begin{center}
\includegraphics[width=0.7\textwidth]{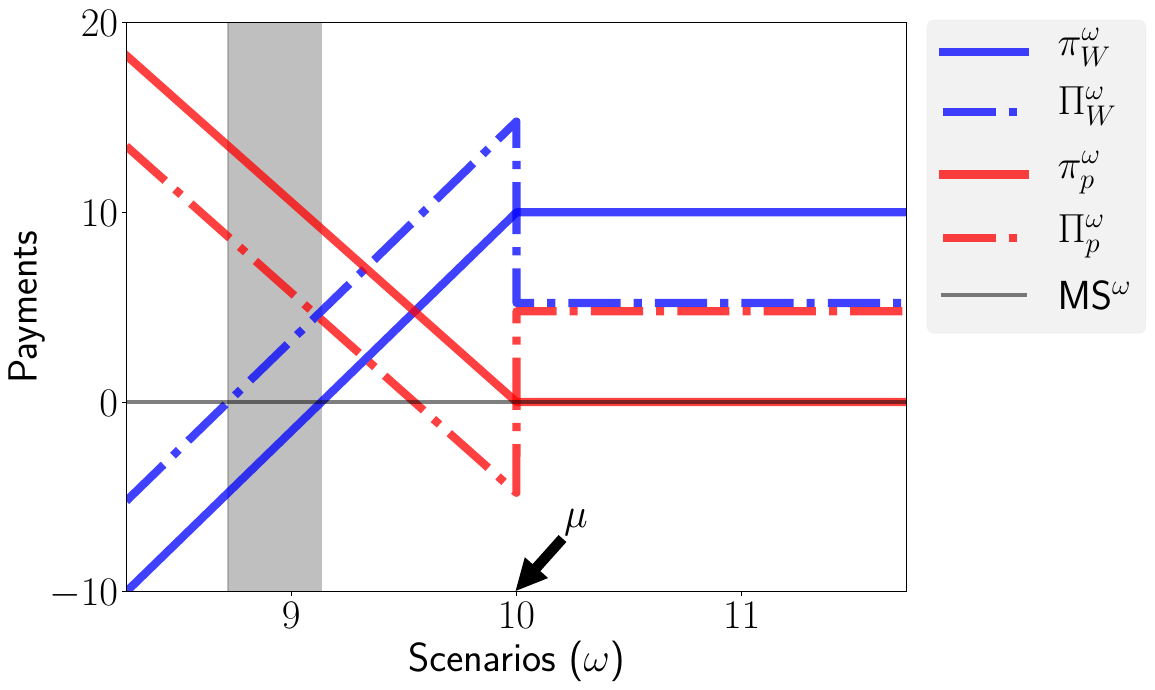}    
\caption{Profits for the copperplate power system example with $\mu=10, \sigma^2=1$, $\Delta^*=\sqrt{3}\sigma$, and $ \rho=\frac{\sqrt{3}}{20}$. The contract allows $W$ to avoid financial loss in the shaded set of scenarios. }
\label{fig:2_Players_Social}
\end{center} 
\end{figure}
Letting $\alpha=0$, the set of acceptable trades for $P$ and $W$ are given by
\begin{equation}
\begin{aligned}
 \Acal_P &= \{( q_P, K_P, \Delta_P)\in\Acal_0:   K_P+ 2q_P\geq 1/\rho \},\\ 
\Acal_W &= \{( q_W, K_W, \Delta_W)\in\Acal_0: K_W+2q_W\leq1/\rho \}.
\end{aligned}
\label{eq:A.PW}
\end{equation}
From the above sets, it is straightforward to infer the feasible set of the insurance contracts in \eqref{eq:DA.opt}, given by $\left( q_W, K_W, \Delta_W \right)=\left( q_P, K_P, \Delta_P \right)=\left( q, K, \Delta \right)$ that satisfies
$$ 2q + K =1/\rho, \quad  \delta^\omega_P=\Delta\ind{\omega \leq \mu}, \quad \Delta \in [0,\sqrt{3}\sigma]. $$
The above trades coincide with the set of all (non-degenerate) Stackelberg equilibria of the bilateral trade between $W$ and $P$ in Proposition \ref{prop:bilat}. Given the objective of the central clearing problem \eqref{eq:DA.opt}, we conclude that  the trade mediated by the market maker finds an equilibrium with the highest aggregate variance reduction. We characterize that reduction in the following result.
\color{black} \begin{proposition}
\label{prop:central}
The optimal solutions of \eqref{eq:DA.opt} for the copperplate power system example are given by
$\left( q^*_W, K^*_W, \Delta^*_W \right)=\left( q^*_P, K^*_P, \Delta^*_P \right)=\left( q^*, K^*, \Delta^* \right)$, where
\begin{gather*}
q^* =\frac{\sqrt{3} \sigma}{4\rho \Delta}-\frac{\sqrt{3} \sigma}{8\Delta}, 
\  
K^*= \frac{1}{\rho}\left(\frac{2\Delta-\sqrt{3}\sigma}{2\Delta}\right) +\frac{\sqrt{3} \sigma}{4 \Delta},
\
\Delta^* \in \left[\frac{\sqrt{3}\sigma (2-\rho)}{4} , \sqrt{3} \sigma \right].
\end{gather*}
Moreover, we have
\begin{align}
    \sum_{i=W,P}\left({\sf var}\left[\Pi^{\omega,*}_i\right]- {\sf var}\left[\pi^\omega_i\right] \right)=-\frac{3\sigma^2}{8}\left(\rho^{-1}-\frac{1}{2}\right)^2 < 0.
    \label{eq:prop4.aggvol}
\end{align}
\end{proposition}

\begin{figure*}
\centering
\includegraphics[width=.8\textwidth]{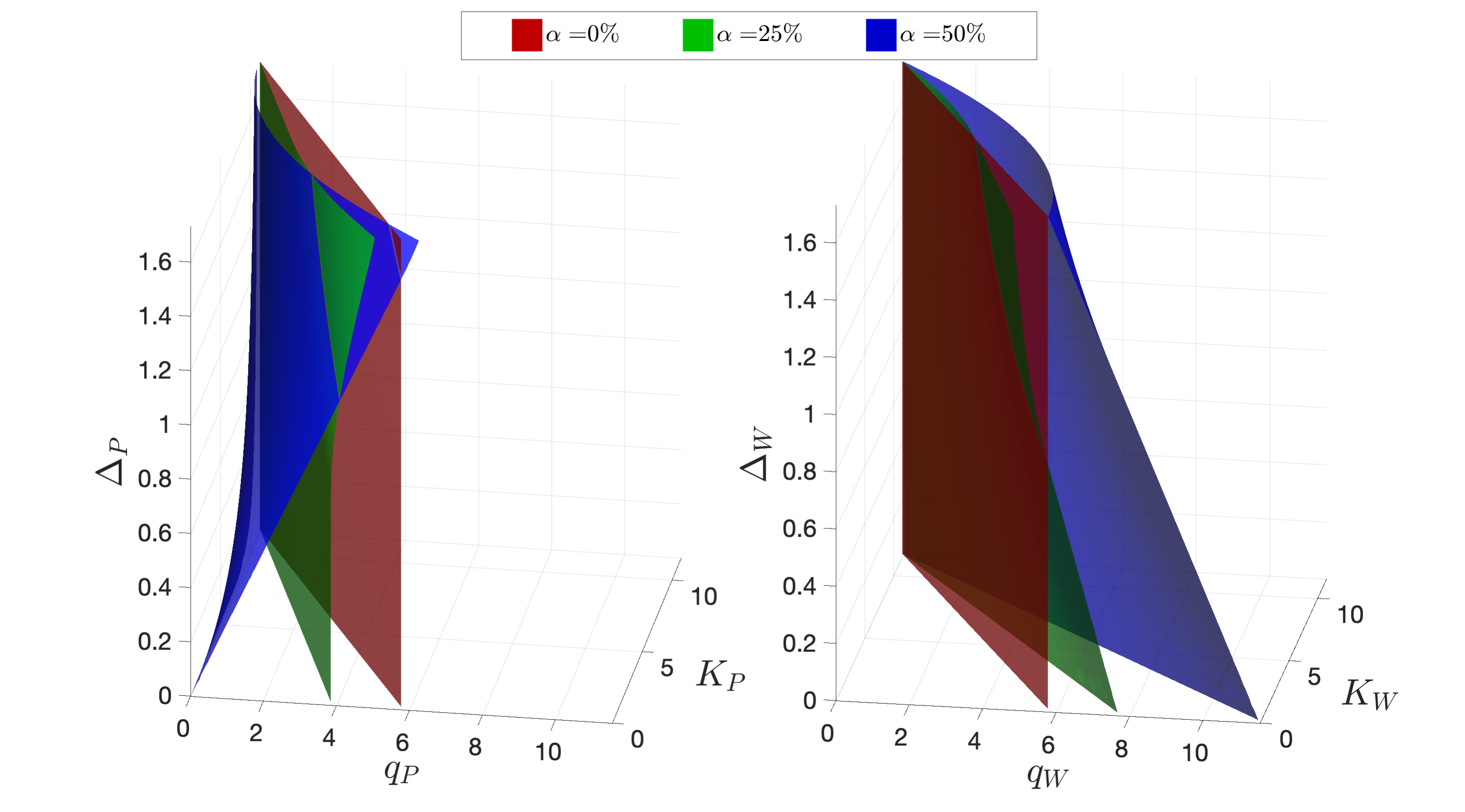}    
\caption{\small{The boundaries of $\Acal_P$ and $\Acal_W$ are portrayed respectively on the left and the right, for the copperplate power system example in Section \ref{sec:central.example}, when $P$ and $W$ both measure risk via ${\sf CVaR}_\alpha$ for different values of $\alpha$. In our experiments, we assume $\mu=10, \sigma^2=1$ and $ \rho=\frac{\sqrt{3}}{20}$, and $\Delta \in [0, \sqrt{3}]$, and compute the sets via the technique outlined in \cite[equation (6)]{Hong}}.}
\label{fig:Aw}
\end{figure*}
{\color{black} Proposition \ref{prop:central} (specifically \eqref{eq:prop4.aggvol}) reveals that aggregate volatilities of $P$ and $W$ strictly decrease as a result of the centralized insurance market. Contrast this result to that in Proposition \ref{prop:bilat}, strict volatility reduction $P$ was only attained for $K>1$. Propositions \ref{prop:bilat} and  \ref{prop:central} reveal that the centralized mechanism picks the subset of Stackelberg equilibria that attain the largest aggregate volatility reduction, while preserving the same expected profits for $P$ and $W$.} Larger $\sigma^2$ implies higher wind uncertainty, leading to largeer variance reduction via our centralized mechanism.
We plot the profits of $W$ and $P$ across the scenarios with the parameters in Figure \ref{fig:2_Players_Social}. Besides decreasing each player's volatility (at no cost to the intermediary), the diagram reveals how $W$ is less exposed to negative profits than without the insurance contract. On the other hand, $P$ is now exposed to negative profits in some scenarios. We will show later in the paper that if $P$ is more risk-averse, she hedges against such losses by requiring higher premium $q$ in the forward stage. 

\subsubsection{The effect of risk-aversion} By varying the risk-aversion parameter, we plot the boundaries of $\Acal_W$ and $\Acal_P$ -- the sets of acceptable trades for the wind power producer and the peaker power plant in our copperplate power system example, respectively -- for various values of $\alpha=\alpha_W=\alpha_P$ in Figure \ref{fig:Aw}.\footnote{The current diagram stands as a correction to \cite[Figure 2]{IFAC}.} Acceptable trades at each $\alpha$ for $W$ lie to the left of the corresponding surface. For $P$, they lie to the right of it. When $\alpha=0$, linearity of expectation allows one to deduce that the acceptability of a trade is independent of $\Delta$. This no longer holds when $\alpha>0$, which is intuitive, as more risk-averse players do not prefer a large $\Delta$. As $\alpha$ grows,  
$P$ requires higher forward premium $q_P$ for a given volume $\Delta_P$. Similar conclusions hold for $W$. She becomes less willing to accept trades with a higher forward premium, the more risk-averse she gets.


\section{Numerical Experiments on the IEEE 14-bus test system}\label{sec:central.example2}

{\color{black}

\begin{figure}
\centering
\includegraphics[width=0.5\textwidth]{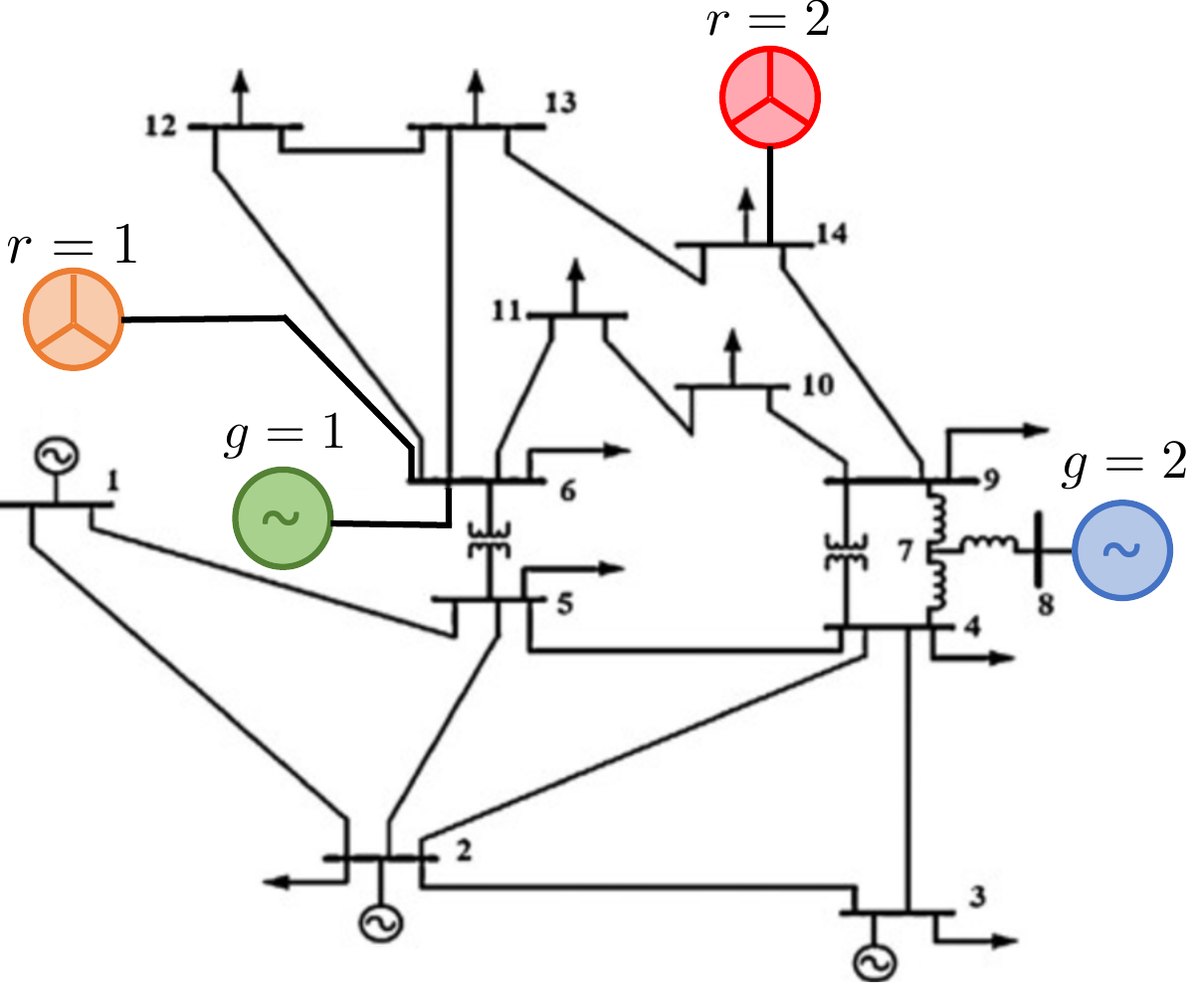}    
\caption{One line diagram of the IEEE 14-bus test system with wind generators added to buses 6 and 14. We consider an insurance market between buyers $r=1,2$ and sellers $g=1,2$.}
\label{fig:IEEE14Bus}
\end{figure}
We now explore the outcomes from the electricity and insurance markets on a modified IEEE 14-bus test system shown in Figure \ref{fig:IEEE14Bus}. Relevant data is adopted from MATPOWER \cite{MATPOWER}.  All transmission lines are assumed to have a capacity of $35$ MW, except that between buses 1 and 2 ($20$ MW) and another one between buses 2 and 4 ($20$ MW). \footnote{The line capacity constraints are added to vary the real-time prices between different buses, which makes the simulation more realistic.} 
Two wind power producers are added to the network at buses 6 and 14. We encode the uncertainty in available wind in $\Omega=[40,60]$ MW and take $\prob$ to be a uniform distribution, i.e., the mean $\mu=50$ MW. We assume zero production costs for the wind generators.
We have the following marginal costs, reflecting offers by producers
\begin{equation}c_g(x)=0.01x^2+40x, \ \ g=1,2. \label{eq:marginal_cost}\end{equation}
We further assume that the true marginal production costs for each generator $g$ is 20 \$/MWh.\footnote{{\color{black}{ We can also assume that the true cost is \eqref{eq:marginal_cost}, which implies truthful bidding. However, in reality, power plants were observed to bid higher costs \cite{Green}. Nevertheless, if true costs are considered, our mechanism remains applicable, as our stylized examples in \cite{IFAC} suggest.}}} Furthermore, we have three real-time prices of interest (one for each bus at which there is a buyer/seller), $p^{\omega,*}_6$, $p^{\omega,*}_8$ and $p^{\omega,*}_{14}$. Each seller/buyer preferences and insurance contract payments are related only to her corresponding bus's day-ahead and real-time prices. Consider an insurance market with the wind power producers at buses $6$ and $14$ as buyers, and the dispatchable generators at buses $6$ and $8$ as sellers. The sellers are generators with higher production costs compared to others in the power system. To deal with uncertainty, we discretize the set $\Omega$ and consider a finite set of scenarios in $\hat{\Omega}=\{40,41,\dots,59,60\}$ MWs.

\subsection{Clearing Procedure and Solution Method}
With our stylized two-stage electricity market described in Section \ref{sec:model}, we have the following profits for each scenario $\omega\in\hat{\Omega}$:
\begin{align*}
    \pi^\omega_{g=1} & \ = P^*_6X^*_{g=1} + p^{\omega,*}_6(x^{\omega,*}_{g=1}-X^*_{g=1})-20x^{\omega,*}_{g=1},\\ 
    \pi^\omega_{g=2} &\ = P^*_7X^*_{g=2} + p^{\omega,*}_7(x^{\omega,*}_{g=2}-X^*_{g=2})-20x^{\omega,*}_{g=2},\\
    \pi^\omega_{r=1} &\ = P^*_6\mu + p^{\omega,*}_6(x^{\omega,*}_{r=1}-\mu), \\
    \pi^\omega_{r=2}  &\ = P^*_{14}\mu + p^{\omega,*}_{14}(x^{\omega,*}_{r=2}-\mu), 
\end{align*}

\noindent where $P^*_n$ denotes the day-ahead price, $X^*_i$ is the day-ahead dispatch for participant $i$, and $x^{\omega,*}_i$ is the real-time dispatch. Recalling the structure of the profit functions after participating in the centralized mechanism, we have: 
\begin{align*}
    \Pi^\omega_{g=1} & \ = \pi^\omega_{g=1}+q_{g=1}\Delta_{g=1}-\left(p^{\omega,*}_6 - K_{g=1}\right)^+\delta^\omega_{g=1},\\ 
    \Pi^\omega_{g=2} &\ = \pi^\omega_{g=2}+q_{g=2}\Delta_{g=2}-\left(p^{\omega,*}_7 - K_{g=2}\right)^+\delta^\omega_{g=2},\\
    \Pi^\omega_{r=1} &\ =  \pi^\omega_{r=1} -q_{r=1}\Delta_{r=1}+\left(p^{\omega,*}_6 - K_{r=1}\right)^+\Delta_{r=1}, \\
    \Pi^\omega_{r=2} & \ =  \pi^\omega_{r=2}-q_{r=2}\Delta_{r=2}+\left(p^{\omega,*}_{14} - K_{r=2}\right)^+\Delta_{r=2}.
    \end{align*}
Equipped with the electricity market outcomes for each $\omega\in\hat{\Omega}$ from MATPOWER, and adopting ${\sf CVaR}$ with risk-aversion parameters $\alpha_i$'s  for all participants, the social market maker solves
$$
\begin{alignedat}{5}
&{\text{minimize}}   \ \ \  {\sf var}[\Pi^\omega_{r=1}]+{\sf var}[\Pi^\omega_{r=2}]+{\sf var}[\Pi^\omega_{g=1}]+{\sf var}[\Pi^\omega_{g=2}], \\
& \text{subject to} \\
& \delta^\omega_{g=1} + \delta^\omega_{g=2} \ = \Delta_{r=1} \ind{p^{\omega,*}_6\geq K_{r=1}}+\Delta_{r=2}\ind{p^{\omega,*}_{14}\geq K_{r=2}},\\
&   \Delta_{g=1}+\Delta_{g=2}  \ = \Delta_{r=1}+\Delta_{r=2},\\
&0 \leq \delta^\omega_{g=1}\leq \Delta_{g=1}, \  \qquad 0 \leq \delta^\omega_{g=2} \leq \Delta_{g=2}, &  \\
& {\sf CVaR}_{\alpha_i} [- \Pi^\omega_i] \leq {\sf CVaR}_{\alpha_i} [ -\pi^\omega_i ], \qquad {\sf MS}^\omega=0, &\\
&\text{for each participant } \ i , \text{ for each scenario }\omega\in\hat{\Omega}.
\end{alignedat} $$

{\color{black}
The decision variables are  $(q_i,  K_i, \Delta_i )$ for each market participant $i$, $\delta^\omega_{g=1}$ and $\delta^\omega_{g=2}$  for each scenario $\omega\in\hat{\Omega}$. The centralized clearing problem is non-smooth and non-convex. We tackle non-smooth functions through smooth surrogates. Precisely, we replace $\ind{x \geq 0}$ and $(x)^+$ by $\left(1+e^{-\beta x} \right)^{-1}$ and $x \left(1+e^{-\beta x}\right)^{-1}$, respectively, with a  large $\beta$. These smooth surrogates slightly relax the problem but they allow us to utilize powerful numerical optimization techniques. Non-convexity makes it challenging to claim convergence of optimization techniques to a global optimum. As evident in the sequel, our simulations demonstrate that solving it using sequential least-squares quadratic programming (SLSQP) returns meaningful results. This is not surprising, given that SLSQP is widely known to work well in practice in solving nonlinear constrained optimization problems \cite{SLSQP2}. The market clearing procedure is implemented as a Jupyter notebook at \cite{tool} that utilizes the SLSQP implementation outlined in \cite{SLSQP}. In our experiments, we use the following values for the set of allowable trades $\Acal_0$ described in Section \ref{sec:centralized}: 
$$\bar{q}=\bar{K}=0.9\max_{\omega}\{p^{\omega,*}_{6},p^{\omega,*}_{7},p^{\omega,*}_{14}\}, \qquad \bar{\Delta}= 10 \ \text{MW}.$$

For the social market maker $\Mcal$, if $\alpha_i=\alpha=0$ for each participant $i$, then, in view of the constraint ${\sf MS}^\omega=0$, 
the constraint $ {\sf CVaR}_{0} [- \Pi^\omega_i] \leq {\sf CVaR}_{0} [ -\pi^\omega_i ]$ can be replaced with $\E [\Pi^\omega_i ] = \E [ \pi^\omega_i ]$ for each $i$. Instead, if $\Mcal$ is profit-maximizing, he would solve the following problem (over the same decision variables).
$$
\begin{alignedat}{5}
&{\text{maximize}}   \ \ \  \E[{\sf MS}^\omega], \\
& \text{subject to} \\
& \delta^\omega_{g=1} + \delta^\omega_{g=2} \ = \Delta_{r=1} \ind{p^{\omega,*}_6\geq K_{r=1}}+\Delta_{r=2}\ind{p^{\omega,*}_{14}\geq K_{r=2}},\\
&   \Delta_{g=1}+\Delta_{g=2}  \ = \Delta_{r=1}+\Delta_{r=2},\\
&0 \leq \delta^\omega_{g=1}\leq \Delta_{g=1}, \  \qquad 0 \leq \delta^\omega_{g=2} \leq \Delta_{g=2}, &  \\
& {\sf CVaR}_{0} [- \Pi^\omega_i] \leq {\sf CVaR}_{0} [ -\pi^\omega_i ],\\
&\text{for each participant } \ i , \text{ for each scenario }\omega\in\hat{\Omega}.
\end{alignedat}$$

\begin{figure}
\centering
\fbox{\includegraphics[width=0.9\textwidth]{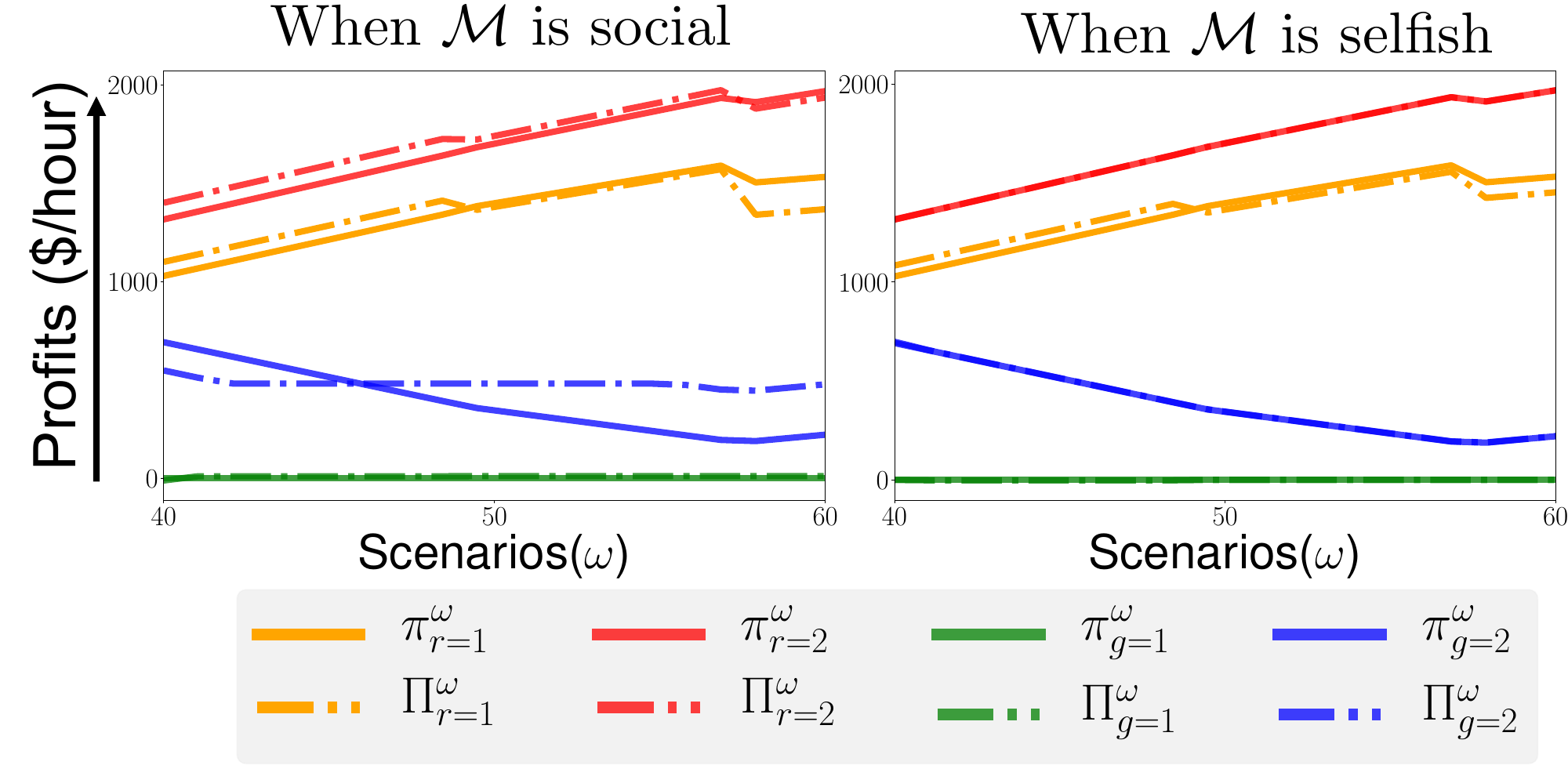}}
\caption{{\color{black}Profits to the buyers/sellers in the IEEE 14-bus system with and without the insurance market. The  figure on the left considers a social intermediary, while the one on the right is derived with a profit-maximizing one. Here, $\alpha=0$ for all participants. }}
\label{fig:IEEEpayments}
\end{figure}

\subsection{Experimental Results}
{\color{black}Figure \ref{fig:IEEEpayments} plots the profits of the market participants with and without the insurance market.
{\color{black} With both social and profit-maximizing $\Mcal$, seller $g=1$ at bus 6 is never dispatched, and hence, receives no profit in the electricity market. Thus, the variance of his profit is zero; it is also kept almost zero after the insurance contract. 
On the contrary, the insurance market reduces the variance for seller $g=2$ drastically (he receives almost the same profit for all scenarios). When $\Mcal$ is a profit-maximizer, there are no guarantees on volatility reduction. Although the market maker here is motivated by maximization of profit, our result in Proposition \ref{prop:profit} showed that a selfish intermediary does not make profits on the average. For the remainder of this section, we focus on the interesting case for which $\Mcal$ is social and report the corresponding optimal numerical insurance values in Table {\ref{table}}. We also report the results when players are risk-averse with $\alpha_i=\alpha=0.5$. When participants are more risk-averse, the quantities cleared are less. For example, $r=1$ does not buy insurance contracts when $\alpha=0.5$. For $\delta^\omega_g$'s, we report a sample of the numerical values because of the large number of scenarios.\footnote{{\color{black}For some $\omega$,  $\sum_g \delta_g^{\omega,*} < \sum_r \Delta_r \ind{p_r^{\omega,*}\geq K_r}$ instead of $\sum_g \delta_g^{\omega,*} = \sum_r \Delta_r \ind{p_r^{\omega,*} \geq K_r}$, due to the smooth approximation of the indicator function. Our results remain largely unaffected as other non-relaxed constraints preserve participants' preferences and consistency.}}  A key feature of our mechanism compared to  bilateral contracts is that insurance quantity $\Delta_r$ bought by buyer $r$ can be matched with multiple sellers through the real-time variables $\delta_g^{\omega}$'s via the market maker's clearing procedure. 
\begin{table}[t]
\centering
\begin{tabular}{cccc|cccc}
    \toprule
     Participant & $q_i $  &  $K_i$ & $\Delta_i$ & $q_i $  &  $K_i$ & $\Delta_i$  \\ 
     $\alpha$ &  0  &  0 &0  & 0.5 & 0.5 &0.5 \\
     \midrule
        $r=1$ & 21.96 & 7.88 & 10 & 23.52 & 7.8 & 0   \\ 
        $r=2$ &17.86 & 16.13 & 10 & 0 & 20.23 & 10  \\ 
        $g=1$ & 0.28 & 36.36 & 10 & 6.78 & 3.64 & 0 \\
        $g=2$ &28.68 & 0 &  10 & 21.73 & 0 &10  \\
        \bottomrule
\end{tabular}
    
\begin{tabular}{cccccc|cccccc}
    \toprule
     $\omega$ & 40 &45  &50  &55   &60  & 40 &45  &50  &55  &60 \\ 
      $\alpha$ & 0 & 0 &  0 & 0& 0& 0.5 & 0.5 &  0.5 & 0.5& 0.5\\
     \midrule
        $\delta^\omega_{g=1}$ & 10 &0.22  &0.22  &3.6  &3.6 &0 &0  &0  &0 &0 \\ 
        $\delta^\omega_{g=2}$ & 10 &7.94  &4.2  &1.5  & 0 & 10 &7.5  &4  &0 &0\\
        \bottomrule
\end{tabular}
        \caption{{\color{black} Optimal insurance contract values for the social market maker. Note that $\delta^\omega_g$'s are proportional to the wind energy shortfall, which is natural. Also, $\Mcal$ varies these values among the sellers depending on the scenario. The more risk-averse the participants become, the less are insurance quantities being cleared. }}
\label{table}
\end{table}

{\color{black}

\subsection{On Complementing Existing Financial Instruments}

Recall that in our mechanism, profit  $\pi^\omega_i$ for market participant $i$ is not limited to the profit from the electricity market alone, but in fact, $\pi^\omega_i$ can represent the payoff from the electricity market in addition to other traditional financial instruments. Next, we demonstrate that our centralized mechanism naturally complements such instruments, leading to significant further volatility reductions while being consistent with participants' preferences. With a slight abuse of notation, assume that buyer $r$ is also a buyer in the bilateral contract described by $(\hat{q},\hat{K},\hat{\Delta})$. Then, he receives a total profit of
$$\Pi^\omega_{r} =  \underbrace{\pi^\omega_{r} -\hat{q}\hat{\Delta}+(\hat{p}^{\omega} - \hat{K})^+\hat{\Delta}}_{\text{electricity market +  bilateral contract}} \  + \quad  \underbrace{- q_r \Delta_r+\left(p^{\omega,*}_b - K_r\right)^+\Delta_r  }_{\text{centralized mechanism}} .$$ 
Similarly, we can define $\Pi^\omega_{g}$. Next, consider the case when there is an existing bilateral cash-settled call option (or insurance contract) between the wind power producer $r=1$ at bus 6 and peaker power plant $g=2$ at bus 7. For illustration, we use the optimal values found for the case when there was no bilateral contract, and assume following bilateral contract values: 
 \begin{gather*}
\hat{q} =  \frac{1}{4} \left(q_{r=1}^*+q_{r=2}^*+q_{g=1}^*+q_{g=2}^*\right), \quad 
\hat{p}^{\omega}=\frac{1}{2}\left({p^{\omega,*}_6+p^{\omega,*}_7}\right), \\
 \hat{K} = \frac{1}{4}\left({K_{r=1}^*+K_{r=2}^*+K_{g=1}^*+K_{g=2}^*}\right), \quad \hat{\Delta} \in \{0, 1, 5, 10\}.
\end{gather*}
Note that when $\hat{\Delta}=0$, the bilateral contract has no effect on the payoffs and we only have the effects of the centralized mechanism's outcomes. In Table \ref{table_bilateral}, we report the percentage change in terms of payoff variances for different values of $\hat{\Delta}$  for all market participants.  Our mechanism complements bilaterally traded call options, leading to further volatility reductions. Also, the bilateral contract did not affect other market participants. 

\begin{table}[t]
    \centering
    \begin{tabular}{ccccc}
    \toprule
     Participant & $\hat{\Delta}=0 $ & $\hat{\Delta}=1$ & $\hat{\Delta}=5$ & $\hat{\Delta}=10$   \\ 
     \midrule
        $r=1$ & -47.7\% & -50.3\% & -56.1\% & -68\% 
        \\ 
        $r=2$ &-28.5\% & -28.5\% & -28.5\% &  -28.5\% \\ 
        $g=1$ & 0\% & 0\% & 0\% & 0\% 
        \\
        $g=2$ &-98.5\% & -99\% & -99.7\% & -100\% 
        \\ 
        \bottomrule
    \end{tabular}
    \caption{{\color{black}The effect of having a bilateral contract  between $r=1$ and $g=2$ in addition to participating in the centralized mechanism. The first column corresponds to the reductions from the centralized mechanism only, and the other columns demonstrate how our mechanism complements bilateral call option trades. Here, $\alpha=0$ for all participants. }}
    \label{table_bilateral}
\end{table}

\begin{figure}
\centering
\includegraphics[width=0.7\textwidth]{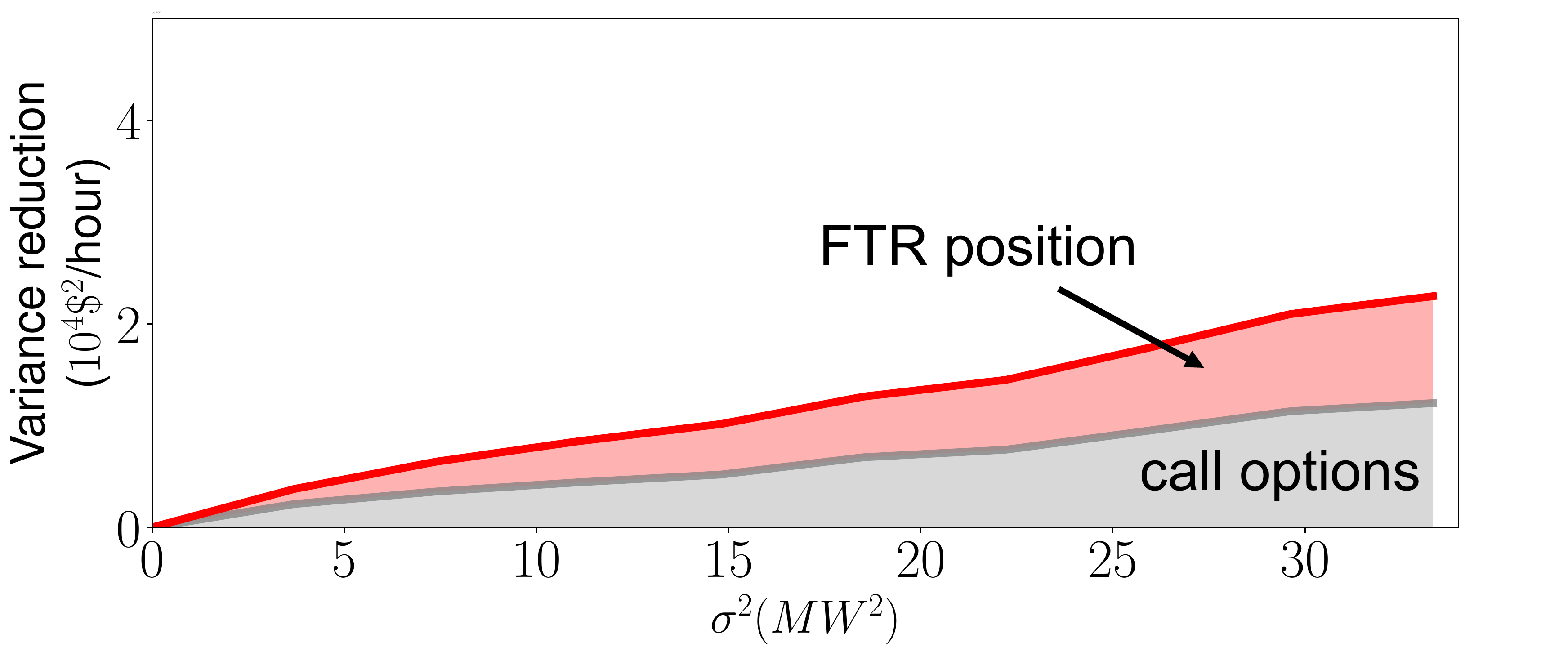}    
\caption{{\color{black}Variance reduction in the profits of  $r=2$ as a function of $\sigma^2$ in the IEEE 14-bus test system. Holding an {\sf FTR} between buses 9 and 14 increases the reduction, due to locational price variation. Here, $\alpha=0$ for all participants.}
\label{fig:IEEE14_Bus_Social_with_FTR}}
\end{figure}

Finally, we also remark that our mechanism does not have any conflicts with other kinds of price variations. For example, a market participant might also hedge against \emph{spatial} price variations using instruments such as financial transmission rights ({\sf FTR}s) \cite{rosellon2013financial}. Holding $f$ {\sf FTR}s between buses $a$ and $b$ entitles a market participant to receive a payment of 
$$ {\sf FTR}^\omega(a,b,f):= (p^{\omega,*}_b-p^{\omega,*}_a) f$$
in scenario $\omega$. Thus, an insurance buyer $r $ located at bus $b$ holding an $\sf FTR$ between buses $a$ and $b$ receives a total profit of
$$\Pi^\omega_{r} =  \underbrace{ \pi^\omega_r+ {\sf FTR}^\omega(a,b,f)}_{\text{electricity market +  spatial hedging}} \  + \quad  \underbrace{-q_r \Delta_r+\left(p^{\omega,*}_b - K_r\right)^+\Delta_r}_{\text{centralized mechanism}} .$$ 
Figure \ref{fig:IEEE14_Bus_Social_with_FTR} illustrates how $r=2$ can reduce its volatility by holding $f =20$ MW worth of {\sf FTR}'s between buses $a=9$ and $b=14$, in addition to the reduction it attains from the insurance contract.}
%

\section{Concluding remarks}
\label{sec:conclusion}
Price volatility in electricity markets is an inevitable consequence of integrating large scale wind energy. In this paper, we have proposed a centralized market for insurance contracts for market participants to tackle the attending financial risks. The centralized mechanism (mediated by a market maker) generalizes bilateral trading of such contracts. On a stylized copperplate power system example, this market provably reduces the profit volatilities of market participants. Numerical experiments on an IEEE 14-bus test system also appear encouraging. For adoption in practice, one needs to estimate the trade volume with real market data from regions with high wind penetration (e.g., Germany, Texas, Denmark). Also, operating such an insurance market in conjunction with current electricity markets will require a carefully designed legal and regulatory framework. {\color{black} Further, our centralized mechanism remains purely financial, but it might motivate market participants to alter their offers/bids in forward electricity markets, such as the day-ahead market. This might change real dispatch schedules, and it would be interesting to formally analyze such effects, which we leave as future research.}

\bibliographystyle{IEEEtran}
\bibliography{references}
\appendix

\section{Proofs}\label{sec:proofs}
\subsubsection*{Proof of Proposition \ref{prop:profit}}
The definitions of $\Acal_g$ and $\Acal_r$ yield
\begin{align*}
\E \[\Pi^\omega_g\] - \E \[ \pi^\omega_g \] \geq 0, \quad
\E \[\Pi^\omega_r \]  - \E \[ \pi^\omega_r \] \geq 0
\end{align*}
for each $g \in \Gfk$ and $r \in \Rfk$. Summing the above inequalities over all $g$ and $r$ yields $\E[{\sf MS}^{\omega}] \leq 0$.
Furthermore, $\E[{\sf MS}^{\omega}] = 0$ is achieved at a feasible point with all $\Delta$'s being identically zero. This completes the proof.
 \subsubsection*{Proof of Proposition \ref{main.2}}
Define $V^\omega_r(q, K, \Delta)
:= \Pi^\omega_r (q, K, \Delta) - \pi^\omega_r$. Then, for each $r\in\Rfk$, we have
 \begin{align*}
&{\sf var}\left[\Pi^\omega_r\right]- {\sf var}\left[\pi^\omega_r\right] =2{\sf cov}(\pi^\omega_r,V^{\omega}_r)+{\sf var}\left[V^{\omega}_r\right]\nonumber \\
& \quad ={\sf cov}(2[P_r^*X^*_r+p^{\omega,*}_r(x^{\omega,*}_r-X^*_r)-c_r(x^{\omega,*}_r)]+V^{\omega}_r,V^{\omega}_r) \\
& \quad ={\sf cov}(2p_r^{\omega,*}(x^{\omega,*}_r-X^*_r) -2c_r(x^{\omega,*}_r) + (p^{\omega,*}_r-K_r)^+\Delta_r, \\
& \quad \qquad \qquad \qquad (p^{\omega,*}_r-K_r)^+\Delta_r)\nonumber \\
& \quad={\sf cov}(2A^\omega_r+B^\omega_r,B^\omega_r).
\end{align*}
The argument for $g \in \Gfk$ is similar and omitted for brevity.

\subsubsection*{Proof of Proposition \ref{prop:bilat}}
Let $P$ choose $\( q, K \) \in \Rset^2_+$. Then, $W$'s payoff from the insurance contract alone is given by
$
V^\omega_W(q, K, \Delta)
:= \Pi^\omega_W (q, K, \Delta) - \pi^\omega_W,$
which, upon utilizing \eqref{eq:Pi}, yields
\begin{align*}
\E\[V^\omega_W(q, K, \Delta)\] 
= \begin{cases}
-q\Delta, & \text{if } K > 1/\rho,\\
- \frac{\Delta}{2} \(  2q +  K - \rho^{-1} \), & \text{otherwise}.
\end{cases}
\end{align*}
We now describe $W$'s best response to $P$'s action.
    \begin{itemize}[leftmargin=*]
    \item If $2q + K < \rho^{-1}$, then $W$ responds by playing $\Delta = \sqrt{3}\sigma$.
    \item If $2q + K = \rho^{-1}$, then $W$ is agnostic to $\Delta$ in $[0, \sqrt{3}\sigma]$.
    \item If $2q + K >\rho^{-1}$, then $W$ chooses $\Delta =  0 $.
    \end{itemize}
Define $V^\omega_P(q, K, \Delta):= \Pi^\omega_P (q, K, \Delta) - \pi^\omega_P,$
as the payoff of $P$ from the insurance contract. Then, the relation in \eqref{eq:Pi} yields
\begin{align}
\E \[ V^\omega_P(q, K, \Delta) \] = - \E \[ V^\omega_W(q, K, \Delta) \].
\label{eq:evp}
\end{align}
Given $W$'s choices, we have the following cases.  
  \begin{itemize}[leftmargin=*]
    \item If $2q + K < \rho^{-1}$, then $\E \[ V^\omega_P(q, K, \Delta) \]<0$. Therefore, $P$ will avoid playing such a $(q, K)$.
    \item If $2q + K = \rho^{-1}$, then $\E \[ V^\omega_P(q, K, \Delta) \]=0$, and $P$ is agnostic to $W$'s choice of $\Delta$ in $[0, \sqrt{3}\sigma]$.
    \item If $2q + K >\rho^{-1}$, then $W$ responds with $\Delta = 0$. And, $P$ receives zero income from the contract.
    \end{itemize}

\noindent Combining them yields the equilibria of $\Gcal$. Now, the difference in variances for $W$ with and without insurance is equal to
 \begin{align}
2{\sf cov}(\pi^\omega_W,V^{\omega,*}_W)+{\sf var}\left[V^{\omega,*}_W\right].  \label{eq:vardiff}\end{align}
When $2q^*+K^*=\rho^{-1}$, we have 
\begin{align*}
& V^{\omega,*}_W(q^*, K^*, \Delta^*) \notag
= \begin{cases}
q^*\Delta^*, & \text{if } \omega \leq \mu,\\
- q^*\Delta^*, & \text{otherwise}.
\end{cases}
\end{align*}
Utilizing $\pi_W^\omega = \mu -( \mu - \omega)^+/\rho$ and $V^{\omega,*}_W$ from the above relation in \eqref{eq:vardiff}, we conclude
\begin{align}&{\sf var}\left[\Pi^\omega_W(q^*, K^*, \Delta^*(q^*, K^*))\right]- {\sf var}\left[\pi^\omega_W\right] \nonumber \\
&\qquad=-(2/\rho){\sf cov}(( \mu - \omega)\ind{\omega<\mu}, V^{\omega,*}_W) + {\sf var}\left[V^{\omega,*}_W\right]\nonumber\\
&\qquad=-\frac{1}{\rho\sqrt{3}\sigma} \int_{\mu-\sqrt{3}\sigma}^{\mu} ( \mu - \omega) q^*\Delta^* d\omega + (q^*\Delta^*)^2\nonumber \\
&\qquad=-\frac{q^*\Delta^*\sqrt{3}\sigma}{2\rho}+ (q^*\Delta^*)^2 \nonumber\\
&\qquad = - q^*\Delta^*\sqrt{3}\sigma \left(q^*+\frac{K^*}{2}\right)+ (q^*\Delta^*)^2\nonumber\\
&\qquad=(q^*)^2\Delta^*(\Delta^*-\sqrt{3}\sigma)-q^*K^*\Delta^*\sqrt{3}\sigma/2.
\label{eq:varineqW}\end{align}
The last expression is nonpositive because $\Delta^*\in[0,\sqrt{3}\sigma]$. 
For $P$, we have $\pi_P^\omega=(\mu-\omega)^+(1/\rho-1)$ and $V^{\omega,*}_P=-V^{\omega,*}_W$. Therefore, similarly, we get
\begin{align}&{\sf var}\left[\Pi^\omega_P(q^*, K^*, \Delta^*(q^*, K^*))\right]- {\sf var}\left[\pi^\omega_P\right] \nonumber \\
&=(q^*)^2\Delta^*(\Delta^*-\sqrt{3}\sigma)-q^*(K^*-1)\Delta^*\sqrt{3}\sigma/2.
\label{eq:varineqP}\end{align}
The rest follows from substituting $\Delta^*=\sqrt{3}\sigma$ in \eqref{eq:varineqW}-\eqref{eq:varineqP}.

\subsubsection*{Proof of Proposition \ref{prop:central}}

The feasible set of \eqref{eq:DA.opt} for the copperplate power system example coincides with the set of nontrivial equilibria of the bilateral trade. We conclude from \eqref{eq:varineqW}-\eqref{eq:varineqP} in the proof of Proposition \ref{prop:bilat} that \eqref{eq:DA.opt} amounts to solving
\begin{equation}
\begin{alignedat}{5}
&{\text{minimize}}   \ \ \ 
	& & 2q^2\Delta(\Delta-\sqrt{3}\sigma)-qK\Delta\sqrt{3}\sigma+q\Delta\sqrt{3}\sigma/2, \\
& \text{subject to}  & & 2q + K = \rho^{-1}, q \geq 0, K \geq 0,\\
&&&  0 \leq \Delta \leq \sqrt{3}\sigma.
\end{alignedat}
\label{eq:DAopt.example}
\end{equation}
Substituting for $K =  \rho^{-1} - 2q$, the objective function of the above problem simplifies to
$ 2q^2 \Delta^2 - q \Delta \sqrt{3} \sigma(\rho^{-1}-1/2)$. 
Being convex quadratic in $q$, it is minimized at 
 $q^*(\Delta) = \min \left\{ \frac{\sqrt{3} \sigma(\rho^{-1}-1/2)}{4 \Delta}, \frac{1}{2\rho} \right\},$
for each $\Delta \in [0, \sqrt{3}\sigma]$. Split the analysis into two cases:
\begin{itemize}[leftmargin=3mm]
\item {Case $\Delta \leq \frac{\sqrt{3}\sigma (2-\rho)}{4} $}: Then, $q^*(\Delta) = \frac{1}{2\rho}$, and the objective function of \eqref{eq:DAopt.example} simplifies to $\frac{1}{2\rho^2} \Delta(\Delta - \sqrt{3} \sigma)+ \frac{\Delta\sqrt{3} \sigma}{4\rho}$. That function is minimized at $\Delta^* = \frac{\sqrt{3}\sigma (2-\rho)}{4} $, taking the value $-\frac{3\sigma^2}{8}(\rho^{-1}-1/2)^2$.

\item {Case $\Delta > \frac{\sqrt{3}\sigma (2-\rho)}{4}$}: Then, we have
$q^*(\Delta) = \frac{\sqrt{3} \sigma(\rho^{-1}-1/2)}{4 \Delta}$ for each $\Delta$, for which the objective function of \eqref{eq:DAopt.example} further simplifies to a constant $-\frac{3\sigma^2}{8}(\rho^{-1}-1/2)^2$.

\end{itemize}
Combining the above two cases and computing the variance reduction at the outcome yields the stated result.

\end{document}